\documentclass{dcds}
 \usepackage{amsmath}
 \usepackage{epsfig} 
 \usepackage{graphics} 

\def\currentvolume{17}
 \def\currentissue{1}
  \def\currentyear{2007}
   \def\currentmonth{January}
    \def\ppages{77--93}

\newtheorem{theo}{Theorem}[section]
\newtheorem{cor}[theo]{Corolary}
\newtheorem{lem}[theo]{Lemma}
\newtheorem{claim}[theo]{Claim}
\newtheorem{prop}[theo]{Proposition}
\theoremstyle{definition}
\newtheorem{defn}[theo]{Definition}
\newtheorem{rem}[theo]{Remark}
\newtheorem{ex}[theo]{Example}

\newcommand{\abs}[1]{\left\vert#1\right\vert}
\newcommand{\set}[1]{\left\{#1\right\}}

\newcommand{\Z}{\mathbb{Z}}
\newcommand{\N}{\mathbb{N}}
\newcommand{\Sg}{\Sigma_\mathbf{A}}
\newcommand{\SgF}{\Sigma_\mathbf{\bar{A}}}
\newcommand{\SgP}{\Sigma_\mathbf{\underline{A}}}
\newcommand{\SgH}{\Sigma_\mathbf{\hat{A}}}
\newcommand{\SgC}{\Sigma_\mathbf{\tilde{A}}}
\newcommand{\s}{\sigma}
\newcommand{\LA}{L_\mathbf{A}}
\newcommand{\LAF}{L_\mathbf{\bar{A}}}
\newcommand{\LAP}{L_{\underline{\mathbf{A}}}}
\newcommand{\LAH}{L_\mathbf{\hat{A}}}
\newcommand{\LAC}{L_\mathbf{\tilde{A}}}
\newcommand{\AZ}{\mathcal{A}^\mathbb{Z}}
\newcommand{\Alf}{\mathcal{A}}
\newcommand{\A}{\mathbf{A}}
\newcommand{\W}{\mathcal{W}}
\newcommand{\F}{\mathcal{F}}
\newcommand{\h}{\mathcal{H}}

\newcommand{\Pp}{\mathcal{P}}

\newcommand{\n}{^{-1}}

\renewcommand{\theenumi}{\roman{enumi}} 

\title[Quasi-Group Shifts]
      {Topological Quasi-Group Shifts}

\author[Marcelo Sobottka]{}


\email{sobottka@dim.uchile.cl}

\thanks{This work was supported by MECESUP UCH0009 and N\'{u}cleo Milenio
Information and Randomness ICM P01-005.}

\begin{document}
\maketitle

\centerline{\scshape Marcelo Sobottka}
\medskip
{\footnotesize
 \centerline{Centro de Modelamiento Matem\'{a}tico}
  \centerline{Facultad de Ciencias F\'{\i}sicas y Matem\'{a}ticas}
   \centerline{Universidad de Chile}
   \centerline{Casilla 170/3-Correo 3, Santiago, Chile}}

\medskip

\begin{abstract}
In this work we characterize those shift spaces which can support
a 1-block quasi-group operation and show the analogous of Kitchens
result: any such shift is conjugated to a product of a full shift
with a finite shift. Moreover, we prove that every expansive
automorphism on a compact zero-dimensional quasi-group that
verifies the medial property, commutativity and has period 2, is
isomorphic to the shift map on a product of a finite quasi-group
with a full shift.
\end{abstract}

\bigskip
\hrule
\noindent
{\footnotesize\em This is a pre-copy-editing, author-produced PDF of an article accepted for publication in Discrete and Continuous Dynamical Systems - Series A (DCDS-A), following peer review. The definitive publisher-authenticated version Marcelo Sobottka, Topological quasi-group shifts. Disc. and Cont. Dynamic. Systems (2007), 17, 1, 77-93, is available online at:\break http://www.aimsciences.org/journals/displayArticles.jsp?paperID=2044 .}
\hrule
\bigskip


\section{Introduction}

One of the main questions concerning symbolic dynamics and
algebraic structures was asked by R. Bowen: characterize group
shifts, that is shifts supporting a group structure so that the
shift map is an automorphism. This question was answered by B.
Kitchens \cite{kitchens}, who showed that any group shift is
conjugated to the product of a full shift with a finite set. A
more general case was studied by N.T. Sindhushayana,  B. Marcus
and M. Trott \cite{marcus}, who proved the analogous result for a
homogeneous shift, that is a shift space $X$ on the alphabet
$\Alf$ for which there exist a group $P(\Alf)$ of permutations of
$\Alf$ and a group shift $Y\subseteq P(\Alf)^\Z$, such that $X$ is
invariant under the action of any element of $Y$.

This work concentrates on quasigroups, often called cancellation
semi-groups, thus with left and right cancellable operations. In
{\S}\ref{II2} we present sufficient and necessary conditions to a
compact zero-dimensional quasi-group $(X,*)$, where is defined an
expansive automorphism $T:X\to X$, to be conjugated and isomorphic
to a Markov shift with a 1-block operation. Furthermore, we give
examples of zero-dimensional quasigroups which verify such
conditions. These are quasi-group versions of results of
\cite{kitchens}, and their proofs use a quasi-group version of
compact zero-dimensional groups (\cite{pont},Theorem 16, pg.77).

We will show that the unique shift spaces which can support a
1-block quasi-group operation are Markov shifts. So, {\S}\ref{II12}
is dedicated to study the case when $\Lambda$ is a Markov shift
and $*$ is a 1-block operation. There, we characterize completely
its structure by supplying a conjugacy with a product of a finite
quasigroup with a full shift.

In the last section we use amalgamations and state splittings
operations (\cite{kitchens},\cite{adler} and \cite{wil}), to
characterize any isomorphism between two quasi-group shifts as in
Kitchens \cite{kitchens}.

\section{Background}

Let $\Alf$ be a finite alphabet and $\AZ$ be the {\it two-sided
full shift} endowed with the product topology (it is a Hausdorff
compact space) . Let $\Lambda\subseteq\AZ$ be a {\it Shift space},
that is a closed shift-invariant set, and denote by
$L_{\Lambda}\subseteq \Alf$ the alphabet used by $\Lambda$.

Let $\W(\Lambda,n)$ be the set of all words or blocks with length
$n$ which are allowed in $\Lambda$ (often we simply write $\W(n)$
instead of $\W(\Lambda,n)$) . Given
$u=[u_1,\ldots,u_n]\in\W(\Lambda,n)$, we write $\F(\Lambda,u)$, or
simply  $\F(u)$, the follower set of $u$:

$$\F(\Lambda,u)=\set{b\in\Alf:[u_1,\ldots,u_n,b]\in\W(\Lambda,n+1)}.$$

In the same way, we define $\Pp(\Lambda,u)$, or simply $\Pp(u)$,
the set of
predecessors of $u$.\\

For $\mathbf{x}=(x_i)_{i\in\Z}\in\Lambda$, $m\leq n$,
we denote $\mathbf{x}[m,n]:=[x_m, x_{m+1},\ldots, x_n]\in\W(n-m+1)$.\\

Let $\s_{\Lambda}$ be the {\it shift map} defined on $\Lambda$,
when the context is clear we
simply put $\s$ instead of $\s_{\Lambda}$.\\

We say that $\Lambda$ is a {\it shift of finite type} (SFT) if
there exists $N\geq 0$ such that for any $\mathbf{x}\in\Lambda$
and for all $n\geq N$ we have
$\F(\mathbf{x}[-n,0])=\F(\mathbf{x}[-N,0])$. In this case we refer
to $\Lambda$ as a $(N+1)$-step SFT.

If $\A$ is a transition matrix on the alphabet $\Alf$, denote by
$\Sg:=\set{\mathbf{x}\in\AZ:\quad A_{x_ix_{i+1}}=1}$\sloppy\  the
{\it two sided Markov shift} and by $\LA$ the alphabet used by
$\Sg$. Without lost of generality, we can assume that all rows and
columns of $\A$ are not null, what is equivalent to say that
$\LA=\Alf$. A Markov shift
is a $1$-step SFT.\\

Let $G$ be a set and $*$ a binary operation on $G$. We say that
$(G,*)$ is a {\it quasigroup} if $*$ is left and right
cancellable:

$$\forall a,b,c\in G,\quad a*b=a*c\quad (or\quad
b*a=c*a)\Longleftrightarrow b=c$$

If, in addition, $G$ is a topological space and $*$ is continuous
with respect to topology of $G$, we say that $(G,*)$ is a {\it
topological quasigroup}. When the context is clear, for $a,b\in
(G,*)$, we write $ab$ instead of $a*b$.

A partition $\mathcal{U}=\set{U_i}_{i\in I}$ of $G$ is said to be
compatible with $*$ if defining $U_i*U_j:=\set{a*b\in G:\quad a\in
U_i,b\in U_j}$\sloppy, so for all $i,j\in I$ there exists $k\in I$
such that $U_i*U_j=U_k$, which is equivalent to say
$(\mathcal{U},*)$ is also a
quasigroup.\\

 Suppose that $(\Lambda,*)$ is a topological quasigroup.
 Then, the shift map is a continuous isomorphism if and only if $*$ is
 given by a $(\ell+r+1)$-block local rule, i.e., there exists
$\ell,r\geq 0$ and $\rho:\W(\ell+r+1)\times\W(\ell+r+1)\to\Alf$,
such that
$$\forall \mathbf{x},\mathbf{y}\in\Lambda,\forall j\in\Z,(x*y)_j
=\rho(\mathbf{x}[j-\ell,j+r],\mathbf{y}[j-\ell,j+r]).$$ In this
case, we say that $\ell$ is the memory and $r$ the anticipation of
$*$. When
$\ell=r=0$, $*$ is a 1-block operation.\\

$(X,T)$ is a {\it topological dynamical system} if $X$ is a
compact space and $T:X\to X$ a homeomorphism. If there exists
$x\in X$, such that $\set{T^n(x):n\geq 0}$ is dense in $X$,
$(X,T)$ is said to be {\it transitive} or {\it irreducible}. Two
topological dynamical systems $(X,T)$ and $(Y,S)$ are {\it
topologically conjugated} if and only if  there exists a
homeomorphism
$\zeta:X\to Y$, such that $\zeta\circ T=S\circ\zeta$.\\

If $(X,*)$ is a topological quasigroup and $(X,T)$ a topological
dynamical system, such that $T:X\to X$ is an automorphism for $*$,
we will denote it by $(X,*,T)$.

We will say that $(X,*,T)$ and $(Y,*,S)$ are isomorphic if and
only if there exists $\zeta:X\to Y$, which is both a topological
conjugation between $(X,T)$ and $(Y,S)$, and an isomorphism
between
$(X,*)$ and $(Y,*)$.\\

If $(X,T)$ is a topological dynamical system, then its topological
entropy \cite{walters} will be denoted by $\mathbf{h}(T)$. When we
refer to the entropy of a shift $(\Lambda,\s_\Lambda)$, we will
write $\mathbf{h}(\Lambda)$.


\section{Expansive automorphisms on zero-dimensional
quasi-groups}\label{II2}

In \cite{kitchens} Kitchens proved that if $(X,\star)$ is a
topological group and $T:X\to X$ is an automorphism, such that:

\begin{description}
  \item[{\sf (H1)}] $X$ is compact (Hausdorff), zero-dimensional and has a numerable topological
  basis, that is, each element $\mathbf{a}\in X$ has a
  clopen fundamental neighborhood $\set{V_n}_{n\geq 1}$:
  $$V_1\supset V_2\supset V_3\supset\cdots\qquad\text{, and}\qquad
  \bigcap_{n=1}^{\infty}V_n=\set{\mathbf{a}}.$$

  \item[{\sf (H2)}] $T$ is an expansive automorphism;
\end{description}

then,

\begin{itemize}
  \item $(X,\star,T)$ is isomorphic by a 1-block code to
  $\bigl(\mathbb{F}\times\Sigma_n,\otimes,\s_\mathbb{F}\times\s_{\Sigma_n}\bigr)$,
  where $\mathbb{F}$ is a finite group with 1-block operation; $\Sigma_n$ is
  a full n shift; and  $\otimes$ is a $k$-block operation, with memory 0 and anticipation
  $k-1$.

  \item If $\mathbf{h}(T)=0$, then $\Sigma_n=\set{a}$, that is,
  the full shift is trivial.

  \item If $T$ is irreducible, then
  $\mathbb{F}=\set{e}$, that is, $\mathbb{F}$ is trivial.

\end{itemize}

Recall that expansivity means that {\it there exists $\mathcal{U}$,
a partition of $X$ by clopen sets (which is finite since $X$ is
compact), such that $\forall \mathbf{x},\mathbf{y}\in X$,
$\mathbf{x}\neq\mathbf{y}$, there exists $n\in\Z$ such that
$T^n(\mathbf{x})$ and $T^n(\mathbf{y})$ belong to distinct sets of
$\mathcal{U}$}.

Our aim is to extend the previous result to the case when
$(X,\star)$ is a topological quasigroup. Now, since a quasigroup has
fewer assumptions about its structure, we need some additional
hypotheses on $(X,\star)$. In particular, it is reasonable to assume
the following property:

\begin{description}
  \item[{\sf (H3)}] $\forall \mathbf{x}\in X$:
  $\mathbf{x}\star X=X\star\mathbf{x}=X$.
\end{description}
(H3) is equivalent to: $\forall \mathbf{y},\mathbf{z}\in X$,
$\exists\mathbf{x_1},\mathbf{x_2}\in X$, such that
$\mathbf{x_1}\star \mathbf{y}=\mathbf{z}$ and
$\mathbf{y}\star\mathbf{x_2}=\mathbf{z}$. Furthermore, since
$\star$ is a quasigroup, the elements $\mathbf{x_1}$ and
$\mathbf{x_2}$ are unique. Notice that if $(X,\star)$ is a finite
quasigroup, then (H3) holds.

Under (H3) we can define on $X$ the following quasi-group
operations $\tilde{\star}$ and $\hat{\star}$:

$$\mathbf{x}\tilde{\star} \mathbf{y}=\mathbf{z}\qquad\Longleftrightarrow\qquad \mathbf{z}\star \mathbf{y}=\mathbf{x}$$

and

$$\mathbf{x}\hat{\star} \mathbf{y}=\mathbf{z}\qquad\Longleftrightarrow\qquad \mathbf{x}\star
\mathbf{z}=\mathbf{y}$$

Also, for any $\mathbf{a}\in X$, we can define the functions
$f_\mathbf{a}:X\to X$ and $f^\mathbf{a}:X\to X$ by
$f_\mathbf{a}(\mathbf{x})= \mathbf{a}\star\mathbf{x}$ and
$f^\mathbf{a}(\mathbf{x})= \mathbf{x}\star\mathbf{a}$. In the same
way we define $\tilde{f}_\mathbf{a}$ and $\tilde{f}^\mathbf{a}$,
using the operation $\tilde{\star}$, and the functions
$\hat{f}_\mathbf{a}$ and $\hat{f}^\mathbf{a}$, using the operation
$\hat{\star}$. It is easy to check that all of these functions are
homeomorphisms.

We recall the identity element plays a fundamental role in the study
of zero-dimensional groups (see \cite{kitchens}, and
\cite{pont},Theorem 16, pg.77). In the case of zero-dimensional
quasi-groups we will need the hypothesis (H3) to define a
substitutive notion:

{\defn\label{II3} For an arbritarily fixed element $\mathbf{e}\in
X$, given $\mathbf{a}\in X$ we define $\mathbf{a}^{-}$ and
$\mathbf{a}^{+}$ as the unique elements in $X$ (which there exist
due (H3)), such that $\mathbf{a}^{-}\star \mathbf{a}=\mathbf{e}$ and
$\mathbf{a}\star \mathbf{a}^{+}=\mathbf{e}$. We say $\mathbf{a}^{-}$
and $\mathbf{a}^{+}$ are respectively the left and right inverses of
$\mathbf{a}$ with respect to $\mathbf{e}$.}

Notice that
$\mathbf{a}^{-}={(f^\mathbf{a})}\n(\mathbf{e})=\tilde{f}_\mathbf{e}(\mathbf{a})$
and
$\mathbf{a}^{+}={f_\mathbf{a}}\n(\mathbf{e})=\hat{f}^\mathbf{e}(\mathbf{a})$.
Moreover, we have that
${(\mathbf{a}^{-})}^{+}={(\mathbf{a}^{+})}^{-}=\mathbf{a}$ and the
maps $\mathbf{a}\mapsto \mathbf{a}^{-}$ and $\mathbf{a}\mapsto
\mathbf{a}^{+}$ are homeomorphisms.

In order to reach our goal, some additional hypotheses over
$(X,\star)$ will be needed.


\subsection{Expansive automorphisms}

Assume that hypotheses (H1) and (H2) hold for $(X,\star,T)$. We
shall prove that there exists a quasi-group shift $(\Lambda,*)$ such
that $(X,\star,T)$ is isomorphic to $(\Lambda,*,\s)$.

{\lem\label{II4} Let $(X,\star,T)$ as above. Then $(X,\star,T)$ is
isomorphic to $(\Lambda,*,\s)$ , where $\Lambda$ is a shift and $*$
is a $k$-block operation, for some $k\geq 1$.}

\begin{proof} From expansibility and 0-dimensionality there exists a
partition $\mathcal{U}$ of $X$, a shift space
$\Lambda\subseteq\mathcal{U}^{\Z}$, and a homeomorphism
$\zeta:X\to\Lambda$, which is a topological conjugacy between
$(X,T)$ and $(\Lambda,\s_\Lambda)$.

In $\Lambda$, we define the quasi-group operation $*$, given by:

$$\forall\mathbf{a},\mathbf{b}\in\Lambda,\qquad
\mathbf{a}*\mathbf{b}:=\zeta\bigl(\zeta^{-1}(\mathbf{a})\star\zeta^{-1}(\mathbf{b})\bigr)$$

We have that $(\Lambda,*)$ is isomorphic by $\zeta$ to
$(X,\star)$. Furthermore, since $\s_\Lambda$ is an automorphism
for $*$, then $*$ is $k$-block, for some $k\geq 1$.
 \end{proof}

In particular we will be interested in the case $*$ being a 1-block
operation. From the proof of Lemma \ref{II4} we deduce that $*$ is a
1-block operation if and only if the partition $\mathcal{U}$ is
compatible with $\star$. For instance, if $X$ is a shift space with
a 1-block operation $\star$, then any partition $\mathcal{U}$ of $X$
by cylinders defined by the same coordinates is compatible with
$\star$ (which means $(X,\star,T)$ is isomorphic to
$(\Lambda,*,\s)$, where $*$ is 1-block).

The natural problem consists in finding such partitions compatible
with the operation for any topological quasigroup in which (H1),
(H2) (and additionally (H3)) hold. This problem remain open.
Therefore, we can ask for the kind of quasi-group structures
allowing to obtain analogous results.

Suppose $(X,\star,T)$ is a topological quasigroup, verifying (H1),
(H2), and such that the following properties hold:

\begin{description}
\item [{\sf (h1)}] $(X,\star)$ is commutative (that is
$f_\mathbf{a}=f^\mathbf{a}$);

\item [{\sf (h2)}] $\star$ has period 2, this means, $\forall \mathbf{a}\in X$,
$f_\mathbf{a}$ has period 2;

\item [{\sf (h3)}] The aforementioned element $\mathbf{e}\in X$ has a fundamental neighborhood system
 $(V_n)_{n\geq 1}$, such that

$$\forall n\geq 1,\qquad\mathbf{e}V_n\subseteq V_n$$

\item [{\sf (h4)}] $(X,\star)$ has the medial property:

$$\forall \mathbf{a},\mathbf{b},\mathbf{c},\mathbf{d}\in X,\qquad
(\mathbf{a}\star\mathbf{b})\star (\mathbf{c}\star
\mathbf{d})=(\mathbf{a}\star\mathbf{c})\star (\mathbf{b}\star
\mathbf{d}).$$\\

\end{description}

Notice that [(h1) and (h2)] is equivalent to  [$\star$,
$\tilde{\star}$ and $\hat{\star}$ are identical]. Thus, for all
$\mathbf{a}\in X$: $\mathbf{a}^{-}=\mathbf{a}^{+}$. Moreover, these
two hypotheses imply that (H3) holds and that the hypothesis (h3) is
equivalent to $\forall n\geq 1$, $\mathbf{e}V_n=V_n\mathbf{e}=V_n$.

Furthermore, from (h3), $\mathbf{e}\star\mathbf{e}=\mathbf{e}$,
which implies $\mathbf{e}^-=\mathbf{e}=\mathbf{e}^+$. We notice
$\mathbf{e}$ is {\em not} an identity element, since in general
$\mathbf{e}\star\mathbf{a}\neq\mathbf{a}\neq\mathbf{a}\star\mathbf{e}$.

Hence, by using (h4), we deduce that for all
$\mathbf{a},\mathbf{b}\in X$, the left inverse with respect to
$\mathbf{e}$ verifies:

\begin{equation}\label{II5}
(\mathbf{a}\star\mathbf{b})^{-}=(\mathbf{a}^{-}\star\mathbf{b}^{-})
\end{equation}

{\ex Let $X=\{0,1\}^\Z\times\{0,1\}^\Z$ and let a 2-block
operation $\star$ with local rule $\rho:\W(X,2)\times\W(X,2)\to
L_X$ be defined for $\mathbf{x},\mathbf{y}\in\W(X,2)$,
$\mathbf{x}=[(x_0^1,x_0^2),(x_1^1,x_1^2)]$ and
$\mathbf{y}=[(y_0^1,y_0^2),(y_1^1,y_1^2)]$, by:

$$\rho(\mathbf{x},\mathbf{y})=\left\{
\begin{array}{ll} (x_0^1+y_0^1,x_0^2+y_0^2+1) &  if\
x_0^1=x_1^1,\ y_0^1=y_1^1\\\\
(x_0^1+y_0^1,x_0^2+y_0^2) &  otherwise\end{array}\right.,$$

where $+$ denote the sum $mod\ 2$.

Then, $(X,\star)$ verifies all of the previous hypotheses, but it is
not a group.

}

The next result is a quasi-group version of a construction done on
topological zero-dimensional groups (\cite{pont},Theorem 17, pg.77). \\

{\theo \label{II6} Let $(X,\star)$ be a topological quasigroup, such
that (H1), (h1), (h2), (h3), (h4) hold. Then, given a neighborhood
$U$ of $\mathbf{e}$, there exists a clopen neighborhood $Q\subseteq
U$ of $\mathbf{e}$, such that $\mathcal{Q}:=\set{\mathbf{a}Q:
\mathbf{a}\in X}$ is a finite partition of $X$ compatible with
$\star$.}

\begin{proof}

{\bf \sf step 1:} Let $(V_n)_{n\geq 1}$ be the neighborhood system
over $\mathbf{e}$, given by (h3). We can suppose $V_{n+1}\subseteq
V_n$, $\forall n\geq 1$.

Let $M\subseteq X$ be a compact subset such that $\mathbf{e}\in
M$. We say that $\mathbf{a}\in M$ can be connected to $\mathbf{e}$
over $M$ by a chain of order $n$, if and only if there exists a
sequence $\mathbf{a}_1=\mathbf{e},\  \mathbf{a}_2,\ \ldots,\
\mathbf{a}_k=\mathbf{a}\in M$, such that

$$\mathbf{a}_i^{-} \mathbf{a}_{i+1}\in V_n,\qquad 1\leq i\leq k-1$$

Let $M_n$ be the set of all points in $M$, which can be connected to
$\mathbf{e}$ over $M$ by a  chain of order $n$. It is
straightforward to see that $M_{n+1}\subseteq M_n$. Moreover, every
point in $M_n$ can be connected to $\mathbf{e}$ over $M_n$ by a
chain of order $n$. In fact, if $\mathbf{a}\in M_n$, then there
exist $\mathbf{a}_1=\mathbf{e},\ \mathbf{a}_2,\ \ldots,\
\mathbf{a}_k=\mathbf{a}\in M$ such that $\mathbf{a}_i^{-}
\mathbf{a}_{i+1}\in V_n$, $1\leq i\leq k-1$. Hence, for any
$j=1,\ldots,k$, $\mathbf{a}_j$ can be connected to $\mathbf{e}$ over
$M$ by a chain of order $n$. Therefore, for any $j=1,\ldots,k$,
$\mathbf{a}_j\in M_n$. Thus, $\mathbf{a}$ can be connected to
$\mathbf{e}$ over $M_n$ by a chain of order $n$.\\

{\bf\sf step 2:} We will show that $M_n$ is a relative open set of
$\tau_{M}:=\set{A\cap M:A\subseteq X\ is\ open}$ the induced
topology in $M$.

Given $\mathbf{a}\in M_n$, we search for a relative open set of
$M$, neighborhood of $\mathbf{a}$, that is a subset of $M_n$.

Let $V':=f^{-1}_{\mathbf{a}^{-}}(V_n)\cap M$, which is a relative
open set of $M$, since the continuity of $f_{\mathbf{a}^{-}}$
implies that $f^{-1}_{\mathbf{a}^{-}}(V_n)$ is open. Moreover,
$f_{\mathbf{a}^{-}}(f^{-1}_{\mathbf{a}^{-}}(V_n))=\mathbf{a}^{-}\star
f^{-1}_{\mathbf{a}^{-}}(V_n)=V_n$, which implies that $\mathbf{a}\in
f^{-1}_{\mathbf{a}^{-}}(V_n)$. Now, let $\mathbf{b}\in V'$, and
$\mathbf{a}_1=\mathbf{e},\ \mathbf{a}_2,\ \ldots,\
\mathbf{a}_k=\mathbf{a}\in M$ be a chain of order $n$ connecting
$\mathbf{a}$ to $\mathbf{e}$ over $M$. Since $\mathbf{b}\in V'$, it
is direct to see that $\mathbf{a}_1=\mathbf{e},\ \mathbf{a}_2,\
\ldots,\ \mathbf{a}_k=\mathbf{a},\ \mathbf{a}_{k+1}=\mathbf{b}\in M$
is a chain of order $n$ connecting $\mathbf{b}$ to $\mathbf{e}$ over
$M$, because
$\mathbf{a}_k^-\mathbf{a}_{k+1}=\mathbf{a}^-\mathbf{b}=f_{\mathbf{a}^-}(\mathbf{b})\in
V_n$. Thus, $\mathbf{b}\in M_n$, so $V'\subseteq M_n$.\\

{\bf\sf step 3:} We need to prove that $M_n$ is a closed set of $X$
or equivalently $M_n^c:=X\setminus M_n$ is open. Since
$M_n^c=M\setminus M_n\cup M^c$, and $M^c$ is open, it is sufficient
to show that $M\setminus M_n$ is open.

Let $\mathbf{a}\in M\setminus M_n$, put
$V':=(f^\mathbf{a})^{-1}(V_n)$, which is an open set since
$f^\mathbf{a}$ is continuous. Notice that
$f^\mathbf{a}(V')=V'\star \mathbf{a}=V_n$ and, since
$\mathbf{e}\in V_n$, it implies that $\mathbf{a}^{-}\in V'$. Let
$V'':={(V')}^{+}:=\set{\mathbf{v''}\in
X:\mathbf{v''}=\mathbf{v'}^{+},\ \mathbf{v'}\in V'}$, which is an
open set containing $\mathbf{a}$ because $\mathbf{a}^-\in
V'$,\linebreak and $(\mathbf{a}^-)^+=\mathbf{a}$.

The set $V''$ cannot intersect $M_n$. In fact, if the contrary
$M_n\cap V''\neq\emptyset$ holds, then it is possible to take
$\mathbf{b}\in M_n\cap V''$ such that $\mathbf{a}_1=\mathbf{e},\
\mathbf{a}_2,\ \ldots,\ \mathbf{a}_{k}=\mathbf{b}\in M$ is a chain
of order $n$ which connects $\mathbf{b}$ to $\mathbf{e}$ over $M$.
But $\mathbf{b}\in V''$, so $\mathbf{b}^{-}\in V'$ and then
$\mathbf{b}^{-} \mathbf{a}\in V_n$. Since $\mathbf{a}\in M$,
$\mathbf{a}_1=\mathbf{e},\ \mathbf{a}_2,\ \ldots,\
\mathbf{a}_{k}=\mathbf{b},\ \mathbf{a}_{k+1}=\mathbf{a}\in M$ is a
chain of order $n$ which connects $\mathbf{a}$ to $\mathbf{e}$
over $M$. Hence, $\mathbf{a}\in M_n$, a contradiction with the
assumption $\mathbf{a}\in M\setminus M_n$.

Since $V''\subseteq M\setminus M_n$ is an open neighborhood of
$\mathbf{a}$, we deduce that $M\setminus M_n$ is an open set.\\

{\bf\sf step 4:} Let $M^*:=\bigcap_{n\geq 1} M_n$. Since
$\set{M_n}_{n\geq 1}$ is a collection of closed sets and $\forall
n\geq 1$, $\mathbf{e}\in M_n$, we have that $\mathbf{e}\in M^*$.

Let us show that $M^*=\set{\mathbf{e}}$. It is sufficient to show
that $\forall \mathbf{b}\in M$, $\mathbf{b}\neq\mathbf{e}$, there
exists $t\geq 1$, such that $\mathbf{b}\notin M_t$, which implies
$\mathbf{b}\notin M^*$.

In fact, given $\mathbf{b}\in M$, and since $M$ is a closed set
and $X$ is a zero-dimensional Hausdorff set, we can write $M=A\cup
B$, a disjoint union of closed sets where $\mathbf{e}\in A$ and
$\mathbf{b}\in B$. Since $A$ and $B$ are compacta, we have that
$A^{-}\star B$ is also compact, so it is a closed set.
Furthermore, $\mathbf{e}\notin A^{-}\star B$ because $A$ and $B$
are disjoint. Then, we can take $V_t$ a neighborhood of
$\mathbf{e}$, such that $V_t\cap (A^{-}\star B)=\emptyset$.

We have that $\mathbf{b}\notin M_t$. In fact, if the contrary
$\mathbf{b}\in M_t$ holds, there would be a chain of order $t$
connecting $\mathbf{b}$ to $\mathbf{e}$ over $M$

$$\mathbf{a}_1=\mathbf{e},\  \mathbf{a}_2,\ \ldots,\
\mathbf{a}_k=\mathbf{b}\in M,\qquad\mathbf{a}_i^{-} \mathbf{a}_{i+1}\in V_t,\qquad 1\leq i\leq k-1$$\\

This would imply that there exists $j$ such that $\mathbf{a}_j\in A$
and $\mathbf{a}_{j+1}\in B$, and so $\mathbf{a}_j^{-}
\mathbf{a}_{j+1}\in  V_t\cap (A^{-}\star B)$.\\

{\bf\sf step 5:} Let $U\subseteq X$ be a neighborhood of
$\mathbf{e}$. Since $\star$ is continuous, there exists $V\subseteq
U$, open neighborhood of $\mathbf{e}$, such that $V\star V\subseteq
U$.

We put $M:=\overline{U}$ (the closure of $U$) and since
$M^*=\set{\mathbf{e}}$, there exists $t\geq 1$, such that
$M_t\subseteq V$. In fact, if there did not exist such $t$, then
for each $n\geq 1$ the set $M_n\cap V^c$ would be closed and not
empty. Then, $\bigcap_{n\geq 1} (M_n\cap V^c)$ would be not empty,
what is a contradiction with $M^*=\set{\mathbf{e}}$.

$M_t$ is a relative open set of $M$, so there exists $W$ an open set
of $X$, such that $M_t=M\cap W=\overline{U}\cap W$. Since
$M_t\subseteq V\subseteq U$, we have that $M_t=M_t\cap W\subseteq
V\cap W\subseteq U\cap W\subseteq \overline{U}\cap W=M_t$, that is,
$M_t=V\cap W$ is an intersection of two open sets, so itself is an
open set.

Let us show that $M_t\star M_t=M_t$. By construction, $M_t\star
M_t\subseteq V \star V\subseteq M$. If $\mathbf{c}\in M_t\star
M_t$, then $\mathbf{c}=\mathbf{a} \mathbf{b}$, with $\mathbf{a},
\mathbf{b}\in M_t$, and there exist two chains of order $t$,
$(\mathbf{a}_i)_{1\leq i\leq k}$ and $(\mathbf{b}_j)_{1\leq j\leq
m}$ which connect respectively $\mathbf{a}$ and $\mathbf{b}$ to
$\mathbf{e}$ over $M_t$.

Thus, we can take the chain $(\mathbf{c}_i)_{1\leq i \leq k+m-1}$:

$$\begin{array}{lcl}
\mathbf{c}_i =  \mathbf{a}_i\mathbf{e}\qquad &,&if\ 1 \leq i \leq k\\
\mathbf{c}_i =   \mathbf{a}\mathbf{b}_{i-k+1}\qquad &,&if\ k+1\leq
i \leq k+m-1
\end{array}$$

Using the medial property and $\mathbf{e}V_t=V_t\mathbf{e}=V_t$ we
get that for all $i\in\set{1,\ldots,\ k+m-1}$ follows
$\mathbf{c}_i^{-} \mathbf{c}_{i+1}\in V_t$. In fact,

$$
\mathbf{c}_i^-\mathbf{c}_{i+1}=\left\{\begin{array}{lcl}
(\mathbf{a}_i^-\mathbf{e})(\mathbf{a}_{i+1}\mathbf{e})=
(\mathbf{a}_i^-\mathbf{a}_{i+1})(\mathbf{e}\mathbf{e})\in
V_t\mathbf{e}
&,&if\ 1 \leq i < k\\
(\mathbf{a}^-\mathbf{e})(\mathbf{a}\mathbf{e})=
(\mathbf{a}^-\mathbf{a})(\mathbf{e}\mathbf{e})\in
V_t &,& if\ i=k \\
(\mathbf{a}^-\mathbf{b}_{i-k+1}^-)(\mathbf{a}\mathbf{b}_{i-k+2})=
(\mathbf{a}^-\mathbf{a})(\mathbf{b}_{i-k+1}^-\mathbf{b}_{i-k+2})\in
\mathbf{e}V_t &,& otherwise
\end{array}\right.$$

Then, this is a chain of order $t$ connecting $\mathbf{c}$ to
$\mathbf{e}$ over $M$, and then $\mathbf{c}\in M_t$.

We have proved $M_t\star M_t\subseteq M_t$. Since $\star$,
$\tilde{\star}$ and $\hat{\star}$ are identical, it follows that
\begin{equation}M_t\star M_t= M_t \qquad (\text{which implies }
M_t^-=M_t=M_t^+)\end{equation}\\

{\bf\sf step 6:} Let $Q:=M_t$, we will show that
$\mathcal{Q}:=\set{\mathbf{a}Q:\mathbf{a}\in X}$ is a partition of
$X$ compatible with $\star$. It follows straightforwardly from
hypothesis (H3) that $\mathcal{Q}$ is an open cover of $X$ and, by
medial property,
$\mathbf{a}Q\star\mathbf{b}Q=(\mathbf{a}\star\mathbf{b})Q$. Then we
only need to prove that $\mathcal{Q}$ is a partition of $X$. To do
this, we introduce the following relation over $X$:

\begin{equation}\label{II7}\mathbf{a}\sim\mathbf{b}\qquad\Longleftrightarrow\qquad\mathbf{a}^{-}
\mathbf{b}\in Q\\\end{equation}

We claim that $\sim$ is an equivalence relation. Clearly $\sim$ is
reflexive. To prove that $\sim$ is symmetric and transitive, we
use that

\begin{equation}\label{II8}\mathbf{a}\sim\mathbf{b}\qquad\Longleftrightarrow\qquad
\mathbf{a}Q=\mathbf{b}Q\\\end{equation}

($\Longrightarrow$) In fact, $\mathbf{a}\sim\mathbf{b}$ is
equivalent to saying that $\mathbf{a}^{-}\star\mathbf{b}\in Q$. So,
for all $\mathbf{q}\in Q$, let $\mathbf{x}\in X$ be such that
$\mathbf{a}\star \mathbf{q}=\mathbf{b}\star \mathbf{x}$. By
multiplying this equation by the left by
$(\mathbf{a}^{-}\star\mathbf{q}^{-})$, and by using the medial
property, we get $\mathbf{e} =
(\mathbf{a}^{-}\star\mathbf{b})\star(\mathbf{q}^{-}\star\mathbf{x})$.
Therefore, since $(\mathbf{a}^{-}\star\mathbf{b})\in Q$, we deduce
that $(\mathbf{q}^{-}\star\mathbf{x})\in Q$ and so $\mathbf{x}\in
Q$. Thus, we can conclude that $\mathbf{a}Q\subseteq\mathbf{b}Q$. By
symmetric reasoning $\mathbf{a}Q\supseteq\mathbf{b}Q$; hence
$\mathbf{a}Q=\mathbf{b}Q$.

($\Longleftarrow$) If $\mathbf{a}Q=\mathbf{b}Q$, then for all
$\mathbf{q}_1\in Q$, there exists $\mathbf{q}_2\in Q$, such that
$\mathbf{a}\star\mathbf{q}_1=\mathbf{b}\star\mathbf{q}_2$. Again,
multiplying this equation by the left by
$(\mathbf{a}^{-}\star\mathbf{q}_1^{-})$, and by using the medial
property, we deduce that $\mathbf{a}^{-}\star\mathbf{b}\in Q$.\\

It is not hard to see that
$\mathcal{Q}=\set{\mathbf{a}Q:\mathbf{a}\in
X}=\set{[\mathbf{b}]:\mathbf{b}\in X}$, where
$[\mathbf{b}]:=\set{\mathbf{c}\in X:\mathbf{c}\sim\mathbf{b}}$ is
the equivalence class of $\mathbf{b}$. Then, $\mathcal{Q}$ is a
partition of $X$ into equivalence classes. Moreover, since $X$ is
compact, $\mathcal{Q}$ is finite.  \end{proof}

 Using last theorem, the following result has a proof similar to
the one of Proposition 2 at \cite{kitchens}.

{\prop\label{II9} Let $(X,\star,T)$ be a topological quasigroup,
such that all hypotheses of Theorem \ref{II6} hold, and $T:X\to X$
is an expansive automorphism (that is, (H2) holds). Then
$(X,\star,T)$ is isomorphic to $(\Lambda,*,\s)$ , where $\Lambda$ is
a shift and $*$ is a 1-block operation.}


 {\rem The hypotheses (h1)-(h4) are strongly restrictive. In fact,
for a finite quasigroup $(X,\star)$, D\'{e}nes-Keedwell
(\cite{denes}, Theorem 2.2.2, p.70) showed that the medial property
implies the quasi-group operation comes from a Abelian group
operation, that is, there exist an Abelian group operation $+$ on
$X$, two automorphisms $\eta$ and $\rho$ on $X$, and $c\in X$, such
that $a\star b = \eta(a)+\rho(b)+c$ for all $a,b\in X$. For our case
of zero-dimensional quasigroups, Theorem \ref{II6}, and propositions
\ref{II9} and \ref{II10}, allow us to get an analogous result
whenever there exists an expansive automorphism (or endomorphism) on
$X$.}


\subsection{1-block quasi-group shifts}

Let $\Lambda\subseteq\AZ$ be a shift space. Suppose that
$(\Lambda,*)$ is a quasigroup where $*$ is a 1-block operation. In
particular, since $*$ is $1$-block, $\s$ is an automorphism over
$(\Lambda,*)$.

{\prop \label{II10} If $(\Lambda,*)$ is as above
 and in addition (H3) holds, then:

\begin{enumerate}
  \item There exists an operation $\bullet$ over $L_\Lambda$, such that
  $(L_\Lambda,\bullet)$ is a quasigroup which induces
  $(\Lambda,*)$;

  \item $\forall k\geq 1$, $\forall g,h\in\W(\Lambda,k)$, we
  have $\abs{\F(g)}=\abs{\F(h)}$.
  Furthermore, if $a\in\F(g)$, then $a\F(h)=\F(g\bullet h)=\F(g)\bullet\F(h)$ and $F(h)a=\F(h\bullet
  g)=\F(h)\bullet\F(g)$;

  \item  $\Lambda$ is a SFT. Moreover $(\Lambda,*,\s)$ is
  isomorphic
  to a Markov shift with a 1-block operation.

\end{enumerate}}

\begin{proof}
  {i.} Since $*$ is a 1-block operation, there exists
  $\rho:L_\Lambda\times L_\Lambda\to L_\Lambda$ a local rule of
  $*$. For $a,b\in L_\Lambda$, put $a\bullet b:=\rho(a,b)$.

  Let us show that $(L_\Lambda,\bullet)$ is a quasigroup. Since $L_\Lambda$ is finite,
  this property is equivalent to the fact
  that for all $a,b\in L_\Lambda$ there exist $c,c'\in L_\Lambda$
  which are the unique
  solutions of $c\bullet a=b$ and $a\bullet c'=b$.

  In fact, if we take
  $\mathbf{y},\mathbf{z}\in\Lambda$, such that $y_0=a$,
  $z_0=b$, there exists $\mathbf{x}\in\Lambda$ a unique solution of
  $\mathbf{x}*\mathbf{y}=\mathbf{z}$. So $x_0\bullet y_0=z_0$,
  which means that $c:=x_0$ is solution of $c\bullet a=b$. Since
  $L_\Lambda$ is finite, if we fix $a$, for each $b$ there exists
  a distinct solution $c$. So, we can deduce that
  $\bullet$ is right permutative. Using the same argument we also deduce the
  left permutativity. Then, $(L_\Lambda,\bullet)$ is a
  quasigroup.

  {ii.} The proof of this fact uses similar arguments as in
  the proofs of Proposition \ref{II15} and Claim \ref{II16} after.

  {iii.} Fix $\mathbf{u}=(u_i)_{i\in\Z}\in\Lambda$. We have
  that $L_\Lambda\supseteq\F(\mathbf{u}[0,0])\supseteq\F(\mathbf{u}[-1,0])
  \cdots\supseteq\F(\mathbf{u}[-n,0])\supseteq\F(\mathbf{u}[-n-1,0])\neq\emptyset$. Since
  $L_\Lambda$ is finite, there exists $N$, such that
  $\F(\mathbf{u}[-n,0])=\F(\mathbf{u}[-N,0])$, for all $n\geq N$.

  Furthermore, if $[g_0,\ldots,g_N]\in\W(\Lambda,N+1)$ and
  $[h_1,\ldots,h_k,g_0,\ldots,g_N]\in\W(\Lambda,N+k+1)$, then
  $\F([h_1,\ldots,h_k,g_0,\ldots,g_N])\subseteq \F([g_0,\ldots,g_N])$
  and they are both cosets of $\F(\mathbf{u}[-N,0])$ (by part (ii)).
  Then, $\F([h_1,\ldots,h_k,g_0,\ldots,g_N])= \F([g_0,\ldots,g_N])$ and we deduce
  that $\Lambda$ is a $(N+1)$-step SFT.\\

  To conclude the proof, we define $\Sg$ as the $(N+1)$-block
  presentation of $\Lambda$, and consider the
  $1$-block quasi-group operation induced by $\Lambda$.
 \end{proof}

We notice that there is no evidence about existence of quasigroups
$(\Lambda,*)$ such that $*$ is a $1$-block operation but (H3) does
not hold.

The previous result implies that if (H3) holds, then $(\Lambda,*)$
is a subquasigroup  of $(\AZ,*)$. The following proposition
reproduce the result of Proposition \ref{II10} using the hypothesis
that $(\Lambda,*)$ is a subquasigroup:

{\prop \label{II11} If $(\Lambda,*)$ is a subquasigroup of
$(\AZ,*)$, where $*$ is a $1$-block operation, then there exists
an operation $\bullet$ over $L_\Lambda$, such that
  $(L_\Lambda,\bullet)$ is a quasigroup which induces
  $(\Lambda,*)$.}

\begin{proof} From the fact of $(\AZ,*)$ is quasigroup follows that for
any constant sequences $\mathbf{a},\mathbf{b}\in\AZ$,
$\mathbf{a}=(\ldots,a,a,a,\ldots)$ and
$\mathbf{b}=(\ldots,b,b,b,\ldots)$, there exist unique
$\mathbf{c},\mathbf{c'}\in\AZ$, $\mathbf{c}=(\ldots,c,c,c,\ldots)$
and $\mathbf{c'}=(\ldots,c',c',c',\ldots)$ solutions of the
equations $$\mathbf{a}*\mathbf{c}=\mathbf{b},\qquad
\mathbf{c'}*\mathbf{a}=\mathbf{b}.$$

Hence, denoting the local rule of $*$ as $\bullet$, for any
$a,b\in L_\Lambda$ the equations $a\bullet c=b$ and $c'\bullet
a=b$ also have unique solutions, which implies $(\LA,\bullet)$ is
quasigroup.
 \end{proof}


\section{Topological 1-block quasi-group Markov
shifts}\label{II12}

In this section we consider the case of subquasigroups
$(\Sg,*,\s)\subseteq (\AZ,*,\s)$, where $\Sg$ is a topological
Markov shift, and $*$ is a 1-block quasi-group operation.
According to the previous section, the operation $*$ over $\Sg$ is
canonically induced by a quasi-group operation $\bullet$ over
$\LA$, such that for any $a,b,a',b'\in\LA$,

\begin{equation}\label{II13}\begin{array}{lrr}
&& a\in\F(b),a'\in\F(b')\Longrightarrow(a\bullet a')\in\F(b\bullet b')\\
\text{and}&&\\
&& a\in\Pp(b),a'\in\Pp(b')\Longrightarrow(a\bullet
a')\in\Pp(b\bullet b')
\end{array}\end{equation}


\subsection{Elementary properties}

{\claim \label{II14} Let $K\subseteq\LA$. Then $\forall g\in\LA$,
$\abs{gK}=\abs{Kg}=\abs{K}$.}\\

\begin{proof} It follows from the bipermutativity of $*$.
 \end{proof}

{\prop \label{II15} Let $(\Sg,*)$ be a 1-block quasi-group shift.
Then,

 \begin{enumerate}
  \item $\forall g,h\in\LA$, $\abs{\F(g)}=\abs{\F(h)}$ and
  $\abs{\Pp(g)}=\abs{\Pp(h)}$
  \item If $s\in\F(r)$, $s\in\Pp(t)$, then
  $$s\F(h)=\F(r\bullet h),\ \F(h)s=\F(h\bullet r),$$
  $$s\Pp(h)=\Pp(t\bullet h),\ \Pp(h)s=\Pp(h\bullet t).$$
\end{enumerate}}
\begin{proof} {i.}  Since $\bullet$ is bipermutative, $\exists r\in\LA$,
such that $r\bullet h=g$. Let $s\in\F(r)$, for all $h'\in\F(h)$ we
have $s\bullet h'\in\F(r\bullet h)=\F(g)$. Then,
$$s\F(h)\subseteq\F(g),$$
and from Claim \ref{II14}, we deduce $\abs{F(h)}\leq\abs{F(g)}$.

Now, let $f_r$ be the permutation over $\LA$, defined by
$f_r(a)=r\bullet a$. There exists $k\in\N$, such that $f_r^k(h)=h$,
so $h=f_r^{k-1}(g)$.

Thus, for $[r,s]\in\W(2)$, $\forall g'\in\F(g)$, we  have
$$\underbrace{[r,s]*(\ldots([r,s]*[g,g']))}_{[r,s]
\ appears\ k-1\ times}=[f_r^{k-1}(g),f_s^{k-1}(g')]=[h,
f_s^{k-1}(g')]\in\W(2)$$

Then, $\forall g'\in\F(g)$, we have $f_s^{k-1}(g')\in\F(h)$ and so
$f_s^{k-1}(\F(g))\subseteq\F(h)$. From Claim \ref{II14}, we have
$$\abs{f_s^{k-1}(\F(g))}=\underbrace{\abs{s*(s*(\ldots(s*\F(g))))}}_{s
\ appears\ k-1\ times}=\abs{\F(g)}\leq\abs{\F(h)}$$

We conclude the aimed equality for the follower sets. Using similar
arguments we deduce the similar equality for the predecessor sets.

{ii.} It is straightforward from part i. and fact (\ref{II13}).
\end{proof}

{\claim\label{II16} For any $a,b\in\LA$ we have that
$\F(a)\bullet\F(b)=\F(a\bullet b)$ and
$\Pp(a)\bullet\Pp(b)=\Pp(a\bullet b)$}\\

\begin{proof} $$\F(a)\bullet\F(b)=\bigcup_{a'\in\F(a)}a'\F(b)=_{(1)}
\bigcup_{a'\in\F(a)}\F(a\bullet b)=\F(a\bullet b),$$

where $=_{(1)}$ follows from  part (ii) of Proposition \ref{II15}.

For the predecessor sets, we use the same argument. \end{proof}

{\defn Let $\Sg$ and $(\LA,\bullet)$ be as before and define

\begin{itemize}
  \item
  $\LAF:=\set{\F(a):a\in\LA}$

  \item
  $\LAP:=\set{\Pp(a):a\in\LA}$

\end{itemize}}

Notice that $\LAF$ and $\LAP$ are both covers of $\LA$. On $\LAF$
and $\LAP$ we consider the operation canonically defined from the
operation on $\LA$ which will be also denoted by $\bullet$. The
Claim \ref{II16} guarantees that $\bullet$ is closed in $\LAF$ and
$\LAP$.

{\prop\label{II17} $(\LAF,\bullet)$ and $(\LAP,\bullet)$ are
quasi-groups.}

\begin{proof} We will only show the result for $(\LAF,\bullet)$, because
the case $(\LAP,\bullet)$ is entirely analogous.

Since $\LAF$ is finite, to prove that $(\LAF,\bullet)$ is right and
left cancellable, is equivalent to prove for all $\F_1,\F_2\in\LAF$,
there exist $\F_i,\F_{j}\in\LAF$, that verify $\F_1\bullet\F_i=\F_2$
and $\F_{j}\bullet\F_1=\F_2$.

We have that $\F_1=\F(a)$ and $\F_2=\F(b)$, for some $a,b\in\LA$. By
bipermutativity in $\LA$, there exist $x,x'\in\LA$ such that
$a\bullet x=b$ and $x'\bullet a=b$. Then, $\F_i:=\F(x)$ and
$\F_j:=\F(x')$ are the solutions for above equations. \end{proof}

{\cor\label{II18} The Markov shift $\Sg$ has disjoint follower (and
predecessor) sets, i.e., $\F(a)\cap\F(b)\neq\emptyset$ if and only
if $\F(a)=\F(b)$.}

\begin{proof} Suppose $\F(a)\cap\F(b)\neq\emptyset$. Let
$r\in\F(a)\cap\F(b)$ and $c\in\LA$. We have

$$\F(a)\bullet\F(c)=\F(a\bullet c)=_{(*)}r\F(c)=\F(b\bullet c)=\F(b)\bullet\F(c),$$

where $=_{(*)}$ is by Proposition \ref{II15}(ii).

Since $(\LAF,\bullet)$ is bipermutative we conclude that
$\F(a)=\F(b)$. \end{proof}

{\cor\label{II19} $(\LAF,\bullet)$ and $(\LAP,\bullet)$ are
isomorphic. In particular, for any $a,b\in\LA$
we have $\abs{\F(a)}=\abs{\Pp(b)}$.}\\

\begin{proof} Let $\tau:\LAF\to\LAP$ defined by $\tau(\F_1)=\Pp(b)$,
where $b\in\F_1$ is an arbitrary element. Let us show that $\tau$
is well defined, i.e., it depends not on the choice of $b$. In
fact,

$$b,b'\in\F_1=\F(a)\Longleftrightarrow \exists
a\in\Pp(b)\cap\Pp(b')\Longleftrightarrow
\Pp(b)\cap\Pp(b')\neq\emptyset \Longleftrightarrow_{(*)}
\Pp(b)=\Pp(b'),$$

where ${(*)}$ is by Corollary \ref{II18}.

Also from above expressions it is direct that $\tau$ is
one-to-one. On
the other hand, is easy to see that $\tau$ is onto.\\

Now, given $\F_1,\F_2\in\LAF$, let $b_1\in\F_1$ and $b_2\in\F_2$.
We have that $b_1\bullet b_2\in\F_1\bullet\F_2$ and

$$\tau(\F_1\bullet\F_2)=\Pp(b_1\bullet b_2)=\Pp(b_1)\bullet
\Pp(b_2)=\tau(\F_1)\bullet\tau(\F_2).$$

To conclude, notice that this isomorphism implies that
$\abs{\LAF}=\abs{\LAP}$. Since $\LAF$ and $\LAP$ are both partitions
of $\LA$, each of them containing sets with the same cardinality
(Proposition \ref{II15}), we deduce that $\abs{\F(a)}=\abs{\Pp(b)}$,
$\forall a,b\in\LA$. \end{proof}

 {\ex Suppose that $*$ is a group operation. In this case, if we
 denote $e\in\LA$ as the identity element, we have that
 $\F(e)=\F(e)\bullet\F(e)$ and $\Pp(e)=\Pp(e)\bullet\Pp(e)$ which
 implies that $(\F(e),\bullet)$ and $(\Pp(e),\bullet)$ are
 subgroups of $(\LA,\bullet)$. Moreover, since $\LAF$ and $\LAP$
 are the sets of cosets of these subgroups, and $(\LAF,\bullet)$ and
 $(\LAP,\bullet)$ are also groups, we conclude that $\F(e)$ and
 $\Pp(e)$ are normal subgroups.}

{\defn Given $a\in\LA$, let $\F(r)\in\LAF$ and $\Pp(t)\in\LAP$ be
such that $a\in\F(r)\cap\Pp(t)$. We define $\h_a:=\F(r)\cap\Pp(t)$
and denote $\LAH:=\set{\h_a:a\in\LA}$.}\\

Notice that $\h_a$ is well defined because for each $a\in\LA$
there exists a unique $\F(r)\in\LAF$ and $\Pp(t)\in\LAP$
satisfying $a\in\F(r)$ and $a\in\Pp(t)$. Moreover, we can write
$\LAH=\set{\F(r)\cap\Pp(t):\quad r,t\in\LA}$.

Consider the operation $\bullet$ over $\LAH$ as in $\LAF$ and
$\LAP$. The Claim \ref{II20} give us that $\bullet$ is closed in
$\LAH$.

{\claim\label{II20} $\forall \h_1,\h_2\in\LA$, $(\h_1\bullet
\h_2)\in\LAH$. Moreover, $\h_{a\bullet
b}=\h_a\bullet\h_b=a\bullet\h_b=\h_a\bullet b$ and for all
$a\in\LA$, $\abs{\h_a}\abs{\LAH}=\abs{\LA}$. }

\begin{proof} Suppose $\h_1=\F(r)\cap\Pp(t)$ and $\h_2=\F(s)\cap\Pp(u)$.
Thus,

\begin{equation}\label{II21}\begin{array}{lcl}
 \h_1\bullet\h_2 & = &(\F(r)\cap\Pp(t))\bullet(\F(s)\cap
\Pp(u))=\bigcup_{g\in (\F(r)\cap\Pp(t))}g\bigl(\F(s)\cap
\Pp(u)\bigr)\\\\
 & = & \bigcup_{g\in (\F(r)\cap\Pp(t))}\bigl(g\F(s)\cap
g\Pp(u)\bigr) =_{(1)}\F(r\bullet s)\cap\Pp(t\bullet u),
\end{array}\end{equation}

where $=_{(1)}$ comes from Proposition \ref{II15}(ii).

Since $\h_1\bullet\h_2$ is a non-empty intersection of sets in
$\LAF$ and $\LAP$, we deduce that it lies in $\LAH$.

Moreover, from definition of $\h_a$ and $\h_b$ it follows that
$(a\bullet b)\in(\h_a\bullet \h_b)$. Then $\h_{a\bullet
b}=\h_a\bullet\h_b$. On the other hand, from equation (\ref{II21}),
we get $\h_a\bullet\h_b=a\bullet\h_b=\h_a\bullet b$. These last
equalities implies that any element of $\LAH$ can be written as the
product of any other element of $\LAH$ by some element of $\LA$,
which implies $\abs{\h_a}\abs{\LAH}=\abs{\LA}$ for any $a\in\LA$.
\end{proof}

{\prop\label{II22} $(\LAH,\bullet)$ is a quasigroup.}

\begin{proof} Use the same argument as in Proposition \ref{II17}. \end{proof}

{\cor\label{II23} $\forall \h_1,\h_2\in\LAH$,
$\h_1\cap\h_2\neq\emptyset\Longleftrightarrow \h_1=\h_2$.}

\begin{proof} The relation ($\Longleftarrow$) is obvious. For the other
one $(\Longrightarrow)$, put $\h_1=\F(r)\cap\Pp(t)$ and
$\h_2=\F(s)\cap\Pp(u)$. Notice that $\h_1\cap\h_2\neq\emptyset$
implies $\F(r)\cap\F(s)\neq\emptyset$ and
$\Pp(t)\cap\Pp(u)\neq\emptyset$. Thus, by Corollary \ref{II18}, we
have $\F(r)=\F(s)$ and $\Pp(t)=\Pp(u)$, and the result follows.
\end{proof}


\subsection{Homomorphisms and isomorphisms}

Fix $e\in\LA$ and let $\h:=\h_e=\F(\bar{x})\cap\Pp(\bar{y})$, where
$\bar{x}\in\Pp(e)$ and $\bar{y}\in\F(e)$. Given $a\in\LA$, define
$a^{-}$ as the element in $\LA$ that verifies $a^{-}\bullet a=e$.

{\defn\label{II24} Let $S:\LAH\to\LA$ be an arbitrary section of
$\LAH$, i.e., an arbitrary map such that
$\forall \h_1\in\LAH$, $S(\h_1)\in\h_1$.}\\

Notice that $\forall \h_1\in\LAH$, ${\h_{S(\h_1)}=\h_1}$.

{\claim\label{II25} $\forall a\in\LA$, $(S(\h_a)^{-}\bullet
a)\in\h$.}

\begin{proof} By definition of $S$ we have $\h_{S(\h_a)}=\h_a$. Then,

$$\begin{array}{lcl}\h_{S(\h_a)^{-}\bullet a} &
=_{(1)}&\h_{S(\h_a)^{-}}\bullet \h_a =_{(2)}
\h_{S(\h_a)^{-}}\bullet \h_{S(\h_a)}\\
& =_{(1)}&\h_{S(\h_a)^{-}\bullet S(\h_a)}=_{(3)}
\h_e=\h,\end{array}$$

where $=_{(1)}$ is by Claim \ref{II20}, $=_{(2)}$ follows from
definition of $S$, and $=_{(3)}$ follows from\linebreak definition
of $^{-}$. \end{proof}

{\prop\label{II26} The map $\phi:\LA\to\LAH\times\h$ given
  by $\phi(a)=\bigl(\h_a,S(\h_a)^{-}\bullet a\bigr)$ is a
  bijection. Moreover, $\phi^{-1}:\LAH\times\h\to\LA$ is given
  by $\phi^{-1}(\h_a,h)=g$, where $g\in\LA$ is the unique
  element such that $S(\h_a)^{-}\bullet g=h_a$.
  }

\begin{proof} To check $\phi$ is one-to-one let $a,b\in\LA$, then

  $$\begin{array}{rcl}\phi(a)=\phi(b)& \Longleftrightarrow & \bigl(\h_a,S(\h_a)^{-}\bullet a)=
  (\h_b,S(\h_b)^{-}\bullet b\bigr)\\\\
  & \Longleftrightarrow & \h_a=\h_b\ \ and\ \
  S(\h_a)^{-}\bullet a=S(\h_b)^{-}\bullet b \Longleftrightarrow
  a=b\end{array}$$

   Since $\LA$ and $\LAH\times\h$ are both
    finite sets with the same cardinality, by Claim \ref{II20}, $\phi$ is also onto.

 Moreover, given $(\h_a,h)\in\LAH\times\h$, let $g\in\LA$ be the unique
  element such that $h=S(\h_a)^{-}\bullet g$. We have that

  $$\begin{array}{rcl}\h_{S(\h_a)^{-}}\bullet\h_{S(\h_a)}&=_{(1)}&\h_{S(\h_a)^{-}\bullet
  S(\h_a)}=\h\\
  &=_{(2)}&\h_{h}=\h_{S(\h_a)^{-}\bullet
  g}=_{(1)}\h_{S(\h_a)^{-}}\bullet\h_g,\end{array}$$

where $=_{(1)}$ is by Claim \ref{II20}, and $=_{(2)}$ is because
$h\in\h$.

Hence, by Proposition \ref{II22}, we get $\h_g=\h_{S(\h_a)}=\h_a$.
Then, $\phi(g)=(\h_g,S(\h_g)^{-}\bullet g)=(\h_a,S(\h_a)^{-}\bullet
g)=(\h_a,h)$. \end{proof}

{\defn\label{II27} Define in $\LAH\times\h$ the operation
  $\diamond$, given by
  $$(\h_1,h_1)\diamond(\h_2,h_2):=\phi\bigl[\phi^{-1}(\h_1,h_1)\bullet\phi^{-1}(\h_2,h_2)\bigr]$$
  }

  Notice that alternatively we can write

  $$(\h_1,h_1)\diamond(\h_2,h_2)=\bigl(\h_1\bullet\h_2,
  S(\h_1\bullet\h_2)^{-}\bullet(g_1\bullet g_2)\bigr),$$

  where $g_1,g_2\in\LA$ are the unique elements which verify

  $$S(\h_1)^{-}\bullet g_1=h_1,$$
  $$S(\h_2)^{-}\bullet g_2=h_2.$$

Notice that on the first coordinate $\diamond$ coincides with
$\bullet$ on $\LAH$.

  {\prop\label{II28} We can identify $(\LA,\bullet)$ to
  $(\LAH\times\h,\diamond)$.}

  \begin{proof} It follows straightforward from the definition of $\diamond$ that $\phi$ is
  an isomorphism between $(\LA,\bullet)$\linebreak and
  $(\LAH\times\h,\diamond)$.
  \end{proof}

{\defn\label{II29} Define the Markov Shift $\SgH$ on the alphabet
$\LAH$, given by transitions:

$$\h_0\to\h_1\Longleftrightarrow\h_1\subseteq\F(\h_0)$$}

The transitions in Definition \ref{II29} can be defined by
$\h_1\subseteq\F(a)$ for any $a\in\h_0$. In fact, if
$\h_0=\F(w)\cap\Pp(z)$, then for all $a\in\h_0$,

\begin{equation}\label{II30}\F(\h_0)=\bigcup_{a'\in\h_0=\F(w)\cap\Pp(z)}\F(a')=_{(1)}\F(a),\end{equation}

where $=_{(1)}$ is due to the fact that for every $a'\in\Pp(z)$,
we have $z\in\F(a')$, hence $\F(a')=\F(a)$ because the follower
sets partition $\LA$, by Corollary \ref{II18}.

Now, consider the map $(a_i)_{i\in\Z}\in\Sg\mapsto
 \bigl(\phi(a_i)\bigr)_{i\in\Z}=\bigl(\h_{a_i},S(\h_{a_i})^{-}\bullet a_i\bigr)_{i\in\Z}\in\SgH
 \times\h^\Z$, which is also denoted as $\phi$.

 We shall check that $\phi:\Sg\to\SgH\times\h^\Z$ is well
 defined, i.e., for every $(a_i)_{i\in\Z}\in\Sg$ we have
 $\bigl(\phi(a_i)\bigr)_{i\in\Z}\in\SgH\times\h^\Z$. Since
 $\phi\bigl((a_i)_{i\in\Z}\bigr)=\bigl(\phi(a_i)\bigr)_{i\in\Z}=
 \bigl(\h_{a_i},S(\h_{a_i})^{-}\bullet a_i\bigr)_{i\in\Z}$, and for all $i\in\Z$
 we have $(S(\h_{a_i})^{-}\bullet a_i)\in\h$ by Claim \ref{II25}, it suffices to verify
 $(\h_{a_i})_{i\in\Z}\in\SgH$. This last property is fulfilled
 because, if
 $a_i\in\F(a_{i-1})$ and $a_i\in\Pp(a_{i+1})$, then
 $\h_{a_i}=\F(a_{i-1})\cap\Pp(a_{i+1})\subseteq\F(a_{i-1})=_{(*)}\F(\h_{a_{i-1}})$,
 where $=_{(*)}$ is by equation (\ref{II30}).

{\prop\label{II31} We can identify $(\Sg,*,\s)$
 to $(\SgH\times\h^\Z,\star,\s)$, where $\star$ is the 1-block operation
 induced by $\diamond$.}

 \begin{proof} Let $\phi:\Sg\to\SgH\times\h^\Z$ be the previous map.

  $\phi:\Sg\to\SgH\times\h^\Z$ is one-to-one because
  its local rule is (Proposition \ref{II28}).

  On the other hand if $(\h_i,h_i)_{i\in\Z}\in\SgH\times\h^\Z$, from Proposition \ref{II28} we
  get
  $(\h_i,h_i)_{i\in\Z}=\bigl(\h_{a_i},S(\h_{a_i})^{-}\bullet a_i\bigr)_{i\in\Z}$.
  So, to deduce that $\phi$ is onto and $\phi^{-1}$ is 1-block, it is sufficient to show that
  $(a_i)_{i\in\Z}\in\Sg$. Now, by definition of $\SgH$, we have $\forall i\in\Z$, $\h_{a_i}\subseteq
  \F(a_{i-1})$, and so $a_i\in\F(a_{i-1})$.

  Therefore, $(\Sg,*,\s)$ is isomorphic to $(\SgH\times\h^\Z,\star,\s)$. \end{proof}

{\cor\label{II32} $(\SgH,*)$ is a quasigroup, where $*$ is
  the operation induced by $\bullet$ over $\LAH$.}

  \begin{proof} $(\SgH,*)$ is a quasigroup because it is a factor of
  $(\SgH\times\h^\Z,\star)$, which is itself a quasigroup because by Proposition
  \ref{II31} says it is isomorphic to $(\Sg,*)$. \end{proof}

{\claim\label{II33} The shift $\SgH$ verifies
$\F(\h_{\bar{x}})\cap\Pp(\h_{\bar{y}})=\set{\h}$, for all
$\bar{x}\in\Pp(e)$ and $\bar{y}\in\F(e)$. }.

\begin{proof} Since $[\bar{x},e,\bar{y}]\in\W(\Sg,3)$, we have
$[\h_{\bar{x}},\h,\h_{\bar{y}}]\in\W(\SgH,3)$. Then,
$\h\in\F(\h_{\bar{x}})\cap\Pp(\h_{\bar{y}})$.

If $\h_1\in\F(\h_{\bar{x}})\cap\Pp(\h_{\bar{y}})$ then
$[\h_{\bar{x}},\h_1,\h_{\bar{y}}]\in\W(\SgH,3)$. Let $a\in\h_1$,
so that $\h_1=\h_a$. By definition of $\SgH$, we have
$\h_a\subseteq \F(\bar{x})$, and so $a\in\F(\bar{x})$.

On the other hand, also from definition of $\SgH$, it follows that
$\h_{\bar{y}}\subseteq \F(a)$. Then, $\bar{y}\in\F(a)$, which is
equivalent to $a\in\Pp(\bar{y})$.

We deduce $a\in\h=\F(\bar{x})\cap\Pp(\bar{y})$, and so we conclude
$\h_a=\h$. \end{proof}

{\defn\label{II34} Define the shift $\SgF$ with alphabet $\LAF$,
and whose transitions are given by:
$$\F_1\to\F_2\Longleftrightarrow\exists g\in\F_1,\ such\ that\
\F(g)=\F_2$$

Let $\theta:\LA\to\LAF$ be the map defined by $\theta(a)=\F(a)$. It
is an onto homomorphism from $(\LA,\bullet)$ to $(\LAF,\bullet)$, by
Proposition \ref{II10}(ii).

Let $(\SgF,*)$ be the quasigroup, with the operation
  $*$ on $\SgF$, induced by the operation $\bullet$ on
  $\LAF$.

We also denote by $\theta$ the map $(a_i)_{i\in\Z}\in\Sg\mapsto
\bigl(\F(a_i)\bigr)_{i\in\Z}\in\SgF$. Let us show that this map is
well defined. Let $(a_i)_{i\in\Z}\in\Sg$, then $\forall i\in\Z,
a_i\in\F(a_{i-1})$. Thus, $\F(a_{i-1})\to\F(a_i)$, i.e.,
$\theta\bigl((a_i)_{i\in\Z}\bigr)=\bigl(\F(a_i)\bigr)_{i\in\Z}\in\SgF$.}

{\prop\label{II35} With the notations above:
\begin{enumerate}

  \item $\theta:\Sg\to\SgF$ is a homomorphism
  from $(\Sg,*)$ onto $(\SgF,*)$ ;

  \item If $\h=\set{e}$, then the element
  $g$ appearing in the definition of $\SgF$ is unique. In such case,
  $\theta$ is an isomorphism between $(\Sg,*,\s)$ and $(\SgF,*,\s)$.

\end{enumerate}}

\begin{proof} {i.} $\theta:\Sg\to\SgF$ is a homomorphism  because its local
rule
  is a homomorphism from $(\LA,\bullet)$ to
  $(\LAF,\bullet)$.

  Let us check that $\theta$ is onto. Let $(\F_i)_{i\in\Z}\in\SgF$, and notice
  that from definition of $\SgF$,
  $\forall i\in\Z$, $\exists a_i\in\F_{i-1}$,
  such that $\F(a_i)=\F_i$. Then, $\exists (a_i)_{i\in\Z}\in\Sg$,
  verifying
  $\theta((a_i)_{i\in\Z})=(\F(a_i))_{i\in\Z}=(\F_i)_{i\in\Z}$.

{ii.} Suppose $\h=\set{e}$. Let $\F_1\to\F_2$ and $g_1,g_2\in\F_1$
be such that
    $\F(g_1)=\F(g_2)=\F_2$.
    Let $a\in\LA$ be such that $\F_1=\F(a)$, and let
    $h\in\F_2=\F(g_1)=\F(g_2)$. Then,

    $$[a,g_1,h],[a,g_2,h]\in\W(\Sg,3)$$

    This implies that $g_1,g_2\in\F(a)\cap\Pp(h)=\h_1$. Since
    $\h_1=b\bullet\h$ for some $b\in\LA$, and since $\h$ is
    unitary, we deduce that $g_1=g_2$. In this case, it is trivial
    to see that there exists
    $\theta^{-1}$. \end{proof}

{\rem\label{II36} If $\h=\set{e}$, then $\theta^{-1}$ is a 2-block
code, with memory 1:
$$\forall (\F_i)_{i\in\Z}\in\SgF,\qquad
\theta^{-1}\bigl((\F_i)_{i\in\Z}\bigr)=(g_i)_{i\in\Z},$$

where for all $i\in\Z$, $g_i\in\F_{i-1}$ is the unique element such
that $\F(g_i)=\F_i$.}

The following theorems are the analogous statements for quasi-groups
as those of theorems stated in \cite{kitchens}. From our previous
results on quasi-groups these theorems have similar proof than those
in \cite{kitchens}.

{\theo\label{II37} Let $(\Sg,*)$ be a quasigroup, where $\Sg$ is a
Markov shift and $*$ is a 1-block operation. Then,
\begin{enumerate}
  \item $(\Sg,*,\s)$ is isomorphic by a 1-block code to
  $\bigl(\mathbb{F}\times\Sigma_n,\otimes,\s_\mathbb{F}\times\s_{\Sigma_n}\bigr)$,
  where $\mathbb{F}$ is a finite quasigroup with
  1-block operation; $\Sigma_n$ is a
  full n shift; and $\otimes$ is a $k$-block quasi-group operation, with memory
  k-1 and anticipation $0$.

  \item $\mathbf{h}(\Sg)=0$ if and only if $\Sigma_n=\set{(\ldots,a,a,a,\ldots)}$ (i.e.,
  the full shift is trivial).

  \item $\Sg$ is irreducible and has a constant
  sequence if and only if
  $\mathbb{F}=\set{e}$ (i.e., $\mathbb{F}$ is unitary).

 \end{enumerate}}

{\theo\label{II38} Let $(\Sg,*)$ be an irreducible Markov shift,
such that $*$ is a 1-block quasi-group operation. Let
$\mathbf{h}(\Sg)=\log(N)$, where $N=p_1^{q_1}\cdots p_r^{q_r}$ is
the prime decomposition of $N$. Then $(\Sg,*,\s)$ is isomorphic to
$(\mathbb{F}\times\Sigma_N,\otimes,\s_\mathbb{F}\times\s_{\Sigma_N})$,
where $\Sigma_N$ is the full $N$ shift and $\otimes$ is at most
$(q_1+\cdots +q_r)$-block, with anticipation 0.}

%
%
%
%
%
%

{\prop Let $(\AZ,*)$ be a quasigroup, where $*$ is induced by a
$1$-block operation $\bullet$ on $\Alf$. Let $\Sg\subset\AZ$ be a
topological Markov chain. Define $\theta:\LA\to\LAF$ by
$\theta(a)=\F(a)$, as before.

Then $\Sg$ is closed under $*$ if and only if $\theta:\LA\to\LAF$ is
an onto homomorphism. Furthermore, $(\Sg,*)$ is irreducible
(transitive) if and only if there exists $a\in\LA$ such that
$\F^k(a)=\LA$ for some $k\geq 0$, where $\F^k(a)$ is defined
inductively by $\F^{n+1}(a)=\bigcup_{h\in\F^n(a)} \F(h)$.}

%
%

{\rem\label{II39} We can define the shift $\SgP$, in the same way
than Definition \ref{II34}:

$$\Pp_1\to\Pp_2\Longleftrightarrow\exists g\in\Pp_2,\ such\ that\
\Pp(g)=\Pp_1 .$$

If we consider $\SgP$ instead of $\SgF$ in Proposition \ref{II35},
we obtain analogous results, but $\theta^{-1}:\SgP\to\Sg$ will be
a 2-block code with anticipation 1.

Moreover, in the Theorems \ref{II37} and \ref{II38}, $\otimes$ will
be a $k$-block operation with memory 0 and anticipation $k-1$.}


\section{Amalgamation and state splitting}\label{II40}

Let $(\Sg,*)$ be a Markov shift with a 1-block quasi-group
operation. As before, denote by $\bullet$ the quasi-group operation
on $\LA$ induced by $*$. We define the four elementary isomorphisms
as in \cite{wil}:

{\sf - State splitting by successors:} Given $a\in\LA$, let
$\h\subseteq\F(a)$ be a subset such that $\LA/_{\h}:=\set{g\h:\quad
g\in\LA}$ is a partition of $\LA$ compatible with $\bullet$. Define
$$\LAC:=\set{(g,\h_{h}):\quad \h_{h}\subseteq \F(g)}\subseteq
\LA\times\LA/_{\h},$$ where $\h_h$ denotes the coset of $\LA/_{\h}$
containing $h$.

Consider on $\LAC$ the operation coinciding with $\bullet$ in each
coordinate, and let $\SgC$ be the shift defined by the following
transitions:
$$(g,\h_{h})\rightarrow (g',\h_{h'})\qquad\Longleftrightarrow\qquad
g'\in\h_{h},$$ which is considered with the operation canonically
induced by $\LAC$.

The state splitting is the 2-block code, $\varphi:\Sg\to\SgC$,
defined by:

$$[g,h]\in\W(\Sg,2)\qquad\mapsto\qquad (g,\h_h)\in\LAC.$$

Notice that $\varphi^{-1}$ is a 1-block code given by
$(g,\h_h)\in\LAC\mapsto g\in\LA$.

The state splitting is an isomorphism between $\Sg$ and $\SgC$.

{\sf - State splitting by predecessors:} It is defined as in the
previous case, but using $\Pp(a)$ instead of $F(a)$.

{\sf - Amalgamation by common predecessors and disjoint successors:}
Given $a\in\LA$, let $\h\subseteq\F(a)$ be a subset, such that
$\LA/_{\h}:=\set{g\h:\quad g\in\LA}$ is a partition of $\LA$
compatible with the operation $\bullet$. Moreover, suppose that
$\forall x\in\LA$, we have $\h\cap\Pp(x)$ is either empty or
unitary.

Define $\LAC:=\LA/_{\h}$, where it is considered the operation
induced by $\bullet$. Let $\SgC$ be the shift given by transitions:
$$\h_g\rightarrow\h_{g'}\qquad\Longleftrightarrow\qquad\exists
h\in\h_g:\quad \h_{g'}\subseteq\F(h),$$ where is defined the
operation induced by $\LAC$.

The amalgamation is the 1-block code, $\varphi:\Sg\to\SgC$, given
by
$$g\in\LA\qquad\mapsto\qquad \h_g\in\LAC .$$

Notice  that $\varphi^{-1}$ is a 2-block code given by
$[\h_g,\h_{g'}]\in\W(\SgC,2)\mapsto h\in\LA$, where $h$ is the
unique element belonging to $\h_g$, such that
$\h_{g'}\subseteq\F(h)$.

It is straightforward to see that the amalgamation is an isomorphism
between $\Sg$ and $\SgC$.

{\sf - Amalgamation by common successors and disjoint predecessors:}
It is defined in the same way than the previous case, but changing
the roles of the predecessor and the follower sets.

{\theo\label{II41} Two quasi-group SFTs, each of them with 1-block
quasi-group operation, are isomorphic if and only if it is possible
to go from one to other by a finite sequence of elementary
isomorphisms.}

 \begin{proof} Let $(\Sg,\s,*)$ and $(\SgF,\s,\bar{*})$ be both quasi-group
shifts with 1-block operations. Let $\phi:\Sg\to\SgF$ be an
isomorphism between them.

Without lost of generality we can consider $\Sg$ a Markov shift
and $\phi$ a 1-block code (in fact, we can take the $N$-block
presentation of $\Sg$, with $N$ sufficiently large). Furthermore,
$(\LA,\bullet)$ and $(\LAF,\bar{\bullet})$ are both quasi-groups
which induce, respectively, the operations $*$ and $\bar{*}$.

We have that $\phi:\LA\to\LAF$ is an onto homomorphism between
$(\LA,\bullet)$ and $(\LAF,\bar{\bullet})$ (notice that the local
rule of the code is also denoted by $\phi$).

Define $\LA/_{\phi^{-1}}:= \set{\phi^{-1}(\{\bar{a}\}):\quad
\bar{a}\in\LAF}$, which is  a partition of $\LA$ compatible with
$\bullet$. This \linebreak property also holds when we consider for
$n\geq 1$, $\phi:\W(\Sg,n)\to\W(\SgF,n)$, that is

$$\W(\Sg,n)/_{\phi^{-1}}=\set{\phi^{-1}(\{[\bar{a}_1,\ldots,\bar{a}_n]\}):\quad
[\bar{a}_1,\ldots,\bar{a}_n]\in\W(\SgF,n)}, $$

which is a partition of $\W(\Sg,n)$ compatible with $\bullet$.

Since $\phi^{-1}$ is a $N$-block code, there exists $m$ , $1\leq
m\leq N$, such that given $[\bar{a}_1,
\ldots,\bar{a}_N]\in\W(\SgF,N)$, $\forall
[a_1,\ldots,a_N],[a'_1,\ldots,a'_N]\in\phi^{-1}(\{[\bar{a}_1,
\ldots,\bar{a}_N]\})$ we have $a_m=a'_m$.

Fix $[\bar{x}_1, \ldots,\bar{x}_N]\in\W(\SgF,N)$,
$[x_1,\ldots,x_N]\in\phi^{-1}(\{[\bar{x}_1, \ldots,\bar{x}_N]\})$
and put $\h:=\F(x_m)\cap\phi^{-1}(\bar{x}_{m+1})$. It follows that
$\h\subseteq\F(x_m)$ and $\LA/_{\h}$ is a partition of $\LA$,
compatible with $\bullet$. Furthermore, for every $b\in\LA$ the set
$\h\cap\Pp(b)$ has at most one element. In fact, if there was more
than one element in $\h\cap\Pp(b)$, then we could find two distinct
sequences in $\Sg$ with the same image by $\phi$, what is a
contradiction about injectivity of this map.

Denote by $\varphi:\Sg\to\SgC$ the amalgamation by common
predecessors and disjoint successors, where $\SgC:=\LA/_\h$. We
recall $\varphi$ is a  $1$-block code and $\varphi^{-1}$ is a
$2$-block code. We define $\tilde{\phi}:\SgC\to\SgF$  the
$1$-block map which has local rule (which we will also denote by
$\tilde{\phi}$), given by $\tilde{\phi}(\h_g):=\phi(g')$ for any
$g'\in\h_g$.

\begin{figure}[h]
\centering
\includegraphics[width=0.3\linewidth=1.0]{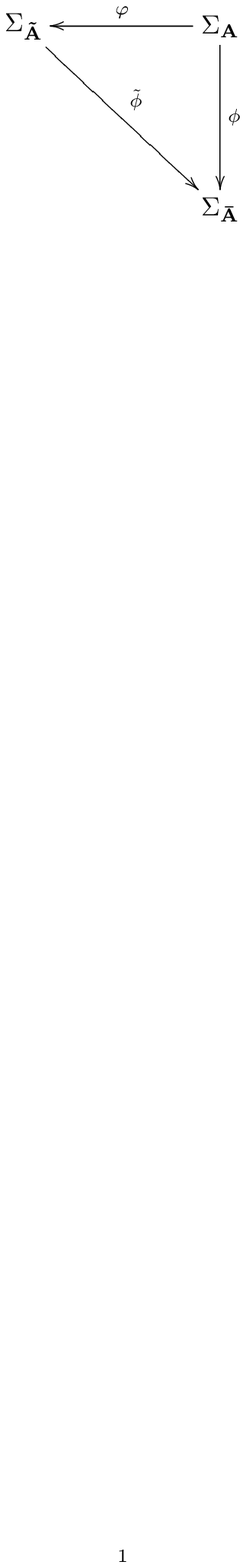}
\end{figure}

Now, we have that given $[\bar{a}_1,
\ldots,\bar{a}_N]\in\W(\SgF,N)$, for all
$[\tilde{a}_1,\ldots,\tilde{a}_N],[\tilde{a}'_1,\ldots,\tilde{a}'_N]
\in\tilde{\phi}^{-1}(\{[\bar{a}_1, \ldots,\bar{a}_N]\})$ follows
that $\tilde{a}_m=\tilde{a}'_m$ and
$\tilde{a}_{m+1}=\tilde{a}'_{m+1}$.

We repeat the above process until we get a $1$-block isomorphism,
$\tilde{\phi}:\SgC\to\SgF$, such that for any $[\bar{a}_1,
\ldots,\bar{a}_N]\in\W(\SgF,N)$, for all
$[\tilde{a}_1,\ldots,\tilde{a}_N],[\tilde{a}'_1,\ldots,\tilde{a}'_N]
\in\tilde{\phi}^{-1}(\{[\bar{a}_1, \ldots,\bar{a}_N]\})$ follows
$\tilde{a}_i=\tilde{a}'_i$, $m\leq i \leq N$.

  To conclude, we come back and applying the amalgamation by common
  successors and disjoint predecessors from the entry $(m-1)$ until
  the first entry. We obtain that $\tilde{\phi}:\SgC\to\SgF$ is a
  1-block isomorphism, and for all $[\bar{a}_1,
  \ldots,\bar{a}_N]\in\W(\SgF,N)$ we have that
  $\tilde{\phi}^{-1}(\{[\bar{a}_1, \ldots,\bar{a}_N]\})$ contains a
  unique $N$-block of $\SgC$. This implies that $\SgC$ and $\SgF$
  are identical.
  \end{proof}

\section*{Acknowledgments}

I want to thank A. Maass, S. Mart\'{\i}nez, and M. Pivato for their
useful discussions and advices on the subject. I also would like to
thank the referee and editor for their comments to strengthen the
presentation of the paper.


\medskip

\medskip

\end{document}